\documentclass[leqno,11pt]{article}
\usepackage{latexsym}
\usepackage{amssymb}
\usepackage{amsfonts}
\usepackage{amsmath}
\usepackage{color}
\usepackage[T1]{fontenc}

\newtheorem{definition}{Definition}[section]
\newtheorem{lemma}[definition]{Lemma}
\newtheorem{theorem}[definition]{Theorem}
\newtheorem{proposition}[definition]{Proposition}
\newtheorem{corollary}[definition]{Corollary}
\newtheorem{remark}[definition]{Remark}
\newtheorem{example}[definition]{Example}

\makeatletter
\long\def\unmarkedfootnote#1{{\long\def\@makefntext##1{##1}\footnotetext{#1}}}
\makeatother

\def\w0{{W_0^{1,A}(\Omega)}}
\def\r{\mathbb R}

\def\ep{\varepsilon}

\def\rn{{{\r}^n}}

\newcommand{\medint}{-\kern  -,395cm\int}
\newcommand{\medintinrigo}{-\kern  -,315cm\int}
\newcommand{\medelle}{-\kern  -,235cm L}
\newcommand{\medellenrigo}{-\kern  -,180cm L}
\newcommand{\qed}{\thinspace\null\nobreak\hfill
\hbox{\vbox{\kern-.2pt\hrule height.2pt
depth.2pt\kern-.2pt\kern-.2pt \hbox to1.8mm {\kern-.2pt\vrule
width.4pt \kern-.2pt\raise1.8mm\vbox to.2pt{} \lower0pt\vtop
to.2pt{}\hfil\kern-.2pt \vrule
width.4pt\kern-.2pt}\kern-.2pt\kern-.2pt \hrule height.2pt
depth.2pt \kern-.2pt}}\par\medbreak}

\title{Regular solutions for nonlinear elliptic equations, 
with convective terms,
in Orlicz spaces}  \numberwithin{equation}{section}

\author{
 Giuseppina Barletta\\
 {\it Dipartimento di Ingegneria Civile, dell'Energia, dell'Ambiente e dei Materiali}
\\ {\it Universit\`a Mediterranea di Reggio Calabria}
 \\ {\it Via Graziella - Loc. Feo di Vito,  89122, Reggio Calabria
 (Italy)}
 \\ {\it e-mail: giuseppina.barletta@unirc.it}
\bigskip
\\
Elisabetta Tornatore
     \\ {\it Dipartimento di Matematica e Informatica }
\\ {\it Universit\`a di Palermo,}
\\  {\it Via Archirafi 34,  90123 Palermo (Italy)} \\
{\it e-mail: elisa.tornatore@unipa.it}
}

\date{}

\begin{document}
	
\maketitle

\begin{abstract} \noindent
We establish some existence and regularity results to the Dirichlet problem, for a class of quasilinear  elliptic equations involving a  partial differential operator, depending on the gradient of the solution.
Our results are formulated in the Orlicz Sobolev spaces and under general growth conditions on the convection term.
The sub and supersolutions method is a key tool in the proof of the existence results.
%We give also some informations on the sign of the solutions we present also some examples in which we highlight the .
\end{abstract}

\unmarkedfootnote {
	\smallskip
\par\noindent {\it
The authors are members of the Gruppo Nazionale  per l'Analisi Matematica, la Probabilit\`{a} e le loro Applicazioni  (GNAMPA) of the Istituto Nazionale di Alta Matematica (INdAM). This research is partially supported by the Ministry of Education, University and Research of Italy, Prin 2017 Nonlinear Differential Problems via Variational, Topological and Set-valued Methods (Project No. 2017AYM8XW).
}}

\section{Introduction}\label{sec1}
Let $\Omega$ be a bounded domain in $\rn$, with $C^{1,\alpha}$ boundary. We consider the following quasilinear elliptic problem involving the $A$-laplacian operator
%\begin{equation}\label{problem}
%\begin{cases}
%- {\rm div}(A'(|\nabla u|)\frac{\nabla u}{|\nabla u|}) =f(x,u, \nabla u) & %{\rm in}\,\,\, \Omega \\
%u  =0 & {\rm on}\,\,\,
%\partial \Omega \,,
%\end{cases}
%\end{equation}
\begin{equation}\label{problem}
\begin{cases}
- \Delta_A u =f(x,u, \nabla u) & {\rm in}\,\,\, \Omega \\
u  =0 & {\rm on}\,\,\,
\partial \Omega \,,
\end{cases}
\end{equation}
where $A:[0,\infty)\to [0,\infty)$ is a convex function, vanishing at $0$, $A\in C^2((0,+\infty))$, 
and $f:\Omega\times\r\times\rn\to\r$ is a Carath\'{e}odory function. The $A$-laplacian operator is defined by $\Delta_A u={\rm div}\left(A'(|\nabla u|)\frac{\nabla u}{|\nabla u|}\right)$. The properties of the function $A$ guarantee that $\Delta_A u$ makes sense also when $|\nabla u|=0$. \\
A wide class of operators can be incorporated in \eqref{problem}. The $p$-Laplacian and the $(p,q)$-Laplacian, for which $A(t)=t^p$ and $A(t)=t^p+t^q$, $t\geq 0$, respectively, are the most known, but we can also consider functions like $A(t)=(\sqrt{1+t^2}-1)^\gamma$, for $t\geq 0$ and $\gamma> 1$ or $A(t)=t^p\lg(1+t)$, for $t\geq 0$ and $p>1$. All the $\Delta_A$ corresponding to the functions $A$ considered above appear in many physical contests, like nonlinear elasticity and plasticity theory.\\ 
The presence of the gradient in the nonlinear term, called convection term, makes variational methods not applicable. Among the techniques used to study problems with a convection term, we cite: topological degree method (\cite{BalFil, Ruiz}), theory of pseudomonotone operators (\cite{GW}), fixed point theorems (\cite{BNV, Zou}), sub and super solution methods (\cite{FMMT, FM, G, NgSch}), approximation methods (\cite{Tanaka}), or a combination of the techniques above (\cite{BaTo1, FMP, MotWin}).\\
We deal with existence, regularity and sign of the solutions to \eqref{problem}. Results in this direction can be found in \cite{BalFil, FMMT, FMP, MotWin, NgSch, Ruiz, Tanaka, Zou}. In the papers above there are various growth conditions on $f$, with respect to each variable, which make it necessary to use different methods to approach the problem, depending on the behavior of the convective term.\\
In all the papers cited above, the abstract framework is the classical Sobolev space $W_0^{1,p}(\Omega)$ and the growth conditions with respect $(s,\xi)\in \r\times \rn$ are of polynomial type. By contrast, we work in Orlicz spaces and take into account a class of operators which, although they depend only on the gradient, cannot be treated in the Sobolev spaces. Furthermore, this allows for $f$ a wider choice than that seen above. Roughly speaking, for a problem with the $p$-Laplacian, a function $f(x,s,\xi)=-c+\frac{|s|^{p^*-1}}{\lg(1+|s|)}+a(|s|)|\xi|^p$ (see Theorem \ref{main2}) is allowed. This does not happen if we consider standard growths.\\
In \cite{BalFil, FMP, Ruiz, Zou} the authors establish the existence and the regularity of positive solutions for a problem involving the $p$-Laplacian. In \cite{BalFil, Ruiz} 
%via a priori bounds and topological degree methods. 
the convection term is a continuous, nonnegative function with subcritical growth with respect to $u$ and growth less then $p$ with respect to $\nabla u$. In \cite{FMP} the growth of the convection term is at most $p-1$ with respect to $u$ and $\nabla u$, while in %and they use z potential estimates and the Shauder fixed point theorem
 \cite{Zou} the convection term is superlinear for $(s,\xi)\to (0,0)$ and its growth is at most $p$ with respect to $u$ and strictly less than $p$ with respect to $\nabla u$. 
%the author prove the existence and the regularity of non negative solutions for a problem with the $p$-laplacian and a .\\ 
An existence and regularity result for the $p$-Laplacian with a convection term that can be singular at $0$ can be found in \cite{NgSch}. The existence of a suitable pair of sub and supersolutions plays a crucial role in their proof. In general, sub and supersolution methods allow to study also the case of a singular convection term, provided the interval of sub and supersolutions does not contain the singular point.\\
%In \cite{BNV} the function $A(t)\sim t^p$ and the convection term is polynomial with respect to $u$ and $\nabla u$.\\
In \cite{MotWin} the authors %prove existence and regularity results, with 
give also sign information on the solutions. They use sub and supersolution methods, in combination with variational techniques, for an operator that can be treated as the $(p,q)$-laplacian. \\ 
Existence and regularity results, for a general operator $A(x,\nabla u)$ can be found in \cite{NgSch, Tanaka}, and in \cite{FMMT} for $A(x,u,\nabla u)$. In \cite{Tanaka} the convection term is a continuous function with growth less than $p-1$ with respect to $u$ and $|\nabla u|$, while in \cite{FMMT, NgSch} the growth is at most $p$ with respect to $|\nabla u|$.\\
%we can find various existence and regularity results obtained via sub and supersolution methods. In particular 
%the authors prove an existence and regularity result for an %operator depending also on $x$, but that can be handled as the $p$-laplacian, 
%of the convection term.\\
% different definitions of sub and supersolutions.
%Among the other paper dealing with \eqref{problem}, we cite also %\cite{Tanaka}, where the the author establishes the existence and the %regularity of a positive solution, via approximation techniques.\\
Let's make some more detailed comments on our new existence and regularity results to \eqref{problem} (Theorems \ref{main}, \ref{main1} and \ref{main2}). 
%involving a general operator $\Delta_A$. $|f|$ is assumed to be bounded from above, with respect to the $\xi \in \rn$, by $A'(|\xi|)|\xi|$.This condition, together with an $L^\infty$ a priori bound for the solution, allows to apply the nonlinear regularity theory (see \cite{Li, Li1} to obtain the regularity of the solutions.\\ In particular, weakening the assumptions on $f$ but strengthening those on $\Omega$, we can prove the existence of a regular solution. Then, under an additional lower bound for $f$ we can apply the strong maximum principle (see \cite{PS}) to affirm that any non negative solution is positive on $\Omega$.\\
In Theorem \ref{main} we assume the existence of an ordered pair of sub and supersolutions $\underline u,\,\overline u \in W^{1,\infty}(\Omega)$ and use Theorem 3.6 in \cite{BaTo1} to prove the existence of a regular solution to \eqref{problem}. The growth condition on $f$, in Theorem \ref{main}, is weaker than that used the Theorem 3.6 in \cite{BaTo1}.
%\textcolor{blue}{Non so se scrivere questo: In this result (as well as in \cite{FMMT, MotWin, NgSch}) the technique used requires only local hypotheses, with respect to the $s$ variable, on the function $f$, related to the pair of  sub and supersolutions.
% which does not explicitly occur because it is a priori bounded.  
%On the other hand....continua sotto}\\
A limit in the use of the method of sub and supersolutions is due to the fact that establishing their existence may not be easy. So we give two existence results, Theorems \ref{main1} and \ref{main2}, where a unified hypothesis on $f$ guarantees the existence of a suitable pair of sub and supersolutions and enable us to apply Theorem \ref{main} to obtain the existence of a regular constant sign solution.\\
The paper is arranged a follows: in Section \ref{sec2} we give the basic definitions and collect some auxiliary results. Our main theorems are proved in Section \ref{sec3}. Finally, in Section \ref{sec4}, we present some examples in which it is easy to verify the existence of constant sub and supersolutions. 

\section{Preliminaries}\label{sec2}

In this Section we give the main definitions on Young functions and define the Orlicz Sobolev spaces that we use in the sequel. For a comprehensive treatment of Young functions and Orlicz spaces we refer the reader to  \cite{Chl, KrRu, RR1, RR2}. We also collect some auxiliary results for the proof of the main theorems.
\begin{definition} A function $A: [0, \infty ) \to [0, \infty]$ is called a Young function if it is convex, vanishes at $0$, and is neither identically equal to $0$, nor to infinity (in $(0,+\infty)$).
\end{definition}
For Young functions
\begin{equation}\label{convexity}
A(\lambda t) \leq \lambda  A(t)\quad  \hbox{for  $\lambda \leq 1$ and $t \geq 0$.}
\end{equation}
\begin{definition}\label{conjygate} The Young conjugate of a Young function $A$ is the Young function $\widetilde A$ defined as
$$\widetilde A (s) = \sup\{ st - A(t):\ t \geq 0\} \quad \hbox{for $s \geq 0$.}$$
\end{definition}
%\begin{definition}
% A Young function $A$ is said to dominate another Young function $B$ near infinity, if there exist constants $c>0$ and $M\geq 0$ such that
%\begin{equation}\label{dyoung}
%   B(t)\leq A(ct)\ \ \hbox{if }\, t\geq M\,.
%\end{equation}
%If \eqref{dyoung} holds with $M=0$, then we say that $A$ dominates $B$ %globally.
 %Two Young functions $A$ and $B$ are called equivalent near infinity %(globally) if they dominate each other near infinity (globally).
%\end{definition}
\begin{definition} A  Young function $A$ is said to satisfy the $\Delta_2$-condition near infinity (briefly $A\in \Delta_2$ near infinity) if it is finite valued and there exist two constants $K\geq 2$ and $M\geq 0$ such that
\begin{equation}\label{delta2young}
A(2t)\leq KA(t)\quad \hbox{for }\ t\geq M\,.
\end{equation}
\end{definition}
\begin{definition} The function $A$ is said to satisfy the $\nabla_2$-condition near infinity (briefly $A\in\nabla_2$ near infinity) if there exist two constants $K>2$ and  $M\geq 0$ such that
\begin{equation}\label{nabla2young}
A(2t)\geq KA(t)\quad \hbox{for }\ t\geq M\,.
\end{equation}
\end{definition}
If \eqref{delta2young} or \eqref{nabla2young} holds with $M=0$, then $A$ is said to satisfy the $\Delta_2$-condition (globally), or the $\nabla_2$-condition (globally), respectively. 
%If $A$ and $B$ are equivalent near infinity then $A\in\Delta_2$ near infinity if and only if $B\in \Delta_2$ near infinity.\\
Given a Young function $A \in C^1([0,+\infty))$, define the quantities 
\begin{equation}\label{i}
	p_A=inf_{t>0}\frac {t\cdot A'(t)}{A(t)},\ \hbox{and}\quad q_A=sup_{t>0}\frac {t\cdot A'(t)}{A(t)}\,.
\end{equation}
The conditions
\begin{equation*}\label{ndg}
	p_A>1\ \hbox{and}\quad q_A<+\infty 
\end{equation*}
are equivalent to the fact that $A\in\nabla_2\cap\Delta_2$ (globally). 

\par
We give basic definitions and the main properties on the Orlicz spaces.
Let $\Omega$ be a measurable set in $\rn$, with $n\geq 1$. Given a Young function $A$, the Orlicz space $L^A(\Omega)$ is the set of all measurable functions $u:\Omega\to\r$ such that the Luxemburg norm
$$\|u\|_{L^A(\Omega)}=\inf\bigg\{\lambda>0:\int_\Omega A\big(\tfrac 1\lambda|u|\big)\,dx\leq 1 \bigg\}$$
is finite. The functional $\|\cdot\|_{L^A(\Omega)}$ is a norm on $L^A(\Omega)$, and the latter is a Banach space (see \cite{Adams}).\\
If $A$ is a Young function, then a generalized  H\"older inequality
\begin{equation}\label{holderyoung}
\int _\Omega |u v|\,dx \leq 2\|u\|_{L^A (\Omega)} \|v\|_{L^{\widetilde A}(\Omega)}
\end{equation}
holds for every $u\in L^A (\Omega)$ and $v\in L^{\widetilde A}(\Omega)$.\\
%Let $A$ and $B$ be  Young functions such that $A$ dominates $B$ near infinity, if %$|\Omega|<\infty$ we have the following embedding
%\begin{equation}\label{embAB}
%   L^A (\Omega)\to L^B (\Omega).
%\end{equation}
%In particular,
%\begin{equation}\label{embA1}
 %  L^A (\Omega)\to L^1 (\Omega).
%\end{equation}
%for any  Young function $A$.
%Also, if $|\Omega|<\infty$,
%\begin{equation}\label{reflex}
%L^A(\Omega)\ \hbox{is reflexive if and only if}\ A\in\Delta_2\cap\nabla_2\ %\hbox{near infinity,} \end{equation}
%see \cite[Corollary 7.2]{Sc}. \\
If $A\in \Delta_2$ globally (or $A\in \Delta_2$ near infinity and $\Omega$ has finite measure) then 
\begin{equation}\label{intfinito}
 \int_{\Omega} A(k|u|)dx <+\infty \ \hbox{for all}\ u\in L^A(\Omega),\ \hbox{all}\ k\geq 0\,. 
\end{equation}
Let $\Omega$ be an open set in $\rn$ with $|\Omega|<\infty$. The isotropic Orlicz-Sobolev spaces 
%$W^{1,A}(\Omega)$ and 
$W^{1,A}_0(\Omega)$ is defined as
%\begin{align*}
%W^{1,A}(\Omega)=\{u:\Omega\to\r:&\, u\hbox{ is weakly differentiable in $\Omega$, and %$u$,\,$|\nabla u| \in L^A (\Omega)$}\}\,
%\end{align*}
%and
\begin{align*}
W^{1,A}_0(\Omega)=\{u:\Omega \to \r: &\, \hbox{the continuation of $u$ by $0$ outside $\Omega$} \\ & \hbox{is weakly differentiable in $\rn$, \, |u|,\, $|\nabla u| \in L^A (\Omega)$}\}.
\end{align*}
The space $W_0^{1,A}(\Omega)$ equipped with the norm 
%
%$W^{1,A}(\Omega)$ equipped with the norm
%$$\|u\|_{W^{1,A}(\Omega)} = \|u\|_{L^A(\Omega)}+\|\nabla u\|_{L^A (\Omega)},$$
%is a Banach space. Also, on the Banach space $W^{1,A}_0(\Omega)$ 
%We can use the equivalent norm 
$$\|u\|_{W^{1,A}_0(\Omega)} = \| |\nabla u| \|_{L^A (\Omega)}.$$
is a Banach space. This norm is equivalent to the standard one
$$\|u\|_{W_0^{1,A}(\Omega)} = \|u\|_{L^A(\Omega)}+\||\nabla u|\|_{L^A (\Omega)}\,.$$
For the Young function $A$ in \eqref{problem}, we assume:
\begin{itemize}
\item [[ A1]] \ $A\in C^2(]0,+\infty[)$ (this  implies $A'\in C^1([0,+\infty[)$);
 \item [[ A2]] there exist two positive constants $\delta$, $g_0>0$ such that
\begin{equation}\label{dg}
    \delta\leq \frac{tA''(t)}{A'(t)}\leq g_0\quad \hbox{for}\ t>0\,.
\end{equation}
\end{itemize}

We point out that, \eqref{dg} guarantees that $A'(0)=0$ and $A\in \nabla_2\cap\Delta_2$ globally. In fact integrating \eqref{dg},
\begin{equation*}
\left(\frac{t}{t_0}\right)^{\delta}\leq \frac{A'(t)}{A'(t_0)}\leq \left(\frac{t}{t_0}\right)^{g_0}\quad \hbox{for}\ t>t_0>0\,.
\end{equation*}
Choosing $t=2t_0$
\begin{equation*}
	2^\delta A'(t_0)\leq A'(2t_0) \leq 2^{g_0}A'(t_0)\quad \hbox{for}\ t>t_0>0\,.
\end{equation*}
Thus
\begin{equation*}
A(2t)=\int_0^{2t}A'(\tau)d\tau=2\int_0^{t}A'(2s)ds\leq 2^{g_0+1}\int_0^{t}A'(s)ds=2^{g_0+1}A(t)\quad \hbox{for all}\; t>0,
\end{equation*}
and
\begin{equation*}
A(2t)=\int_0^{2t}A'(\tau)d\tau=2\int_0^{t}A'(2s)ds\geq 2^{\delta+1} \int_0^{t}A'(s)ds=2^{\delta+1} A(t)\quad \hbox{for all}\; t>0.
\end{equation*}
\par
We investigate the existence and the regularity of the solutions to problem \eqref{problem}. The proof of the existence is based on sub and supersolution methods, while the main tool for the regularity is Theorem 1.7 of \cite{Li1} and the remark immediately after the statement (see also \cite[Theorem 1]{Li}), that we recall below.
\begin{proposition}(see \cite[Theorem 1.7]{Li1})\label{Lib}
Let $\Omega$ be a bounded domain in $\r^n$ with $C^{1,\alpha}$ boundary, for some $0<\alpha\leq 1$. Let $g:[0,+\infty[\to [0,+\infty[$ be a $C^1$, increasing function, satisfying $0<\delta\leq \frac{tg'(t)}{g(t)}\leq g_0$, for $t>0$, and let $G(t)=\int_0^t g(\tau)d\tau$. Consider the problem
$$div({\cal A}(x,u,\nabla u))+B(x,u,\nabla u)=0\ \hbox{in}\ \Omega\,.$$
Suppose ${\cal A}$ and $B$ satisfy the structure conditions (here $a_{ij}(x,z,\eta)=\frac{\partial {\cal A}^i}{\partial\eta_j}$) 
\begin{itemize}
	\item[$(a)$] $\sum_{i,j=1}^na_{ij}(x,z,\eta)\xi_i\xi_j\geq  \frac{g(|\eta|)}{|\eta|}|\xi|^2$
\item[$(b)$] $\sum_{i,j=1}^n|a_{ij}(x,z,\xi)|\leq  \Lambda\frac{g(|\xi|)}{|\xi|}$
\item[$(c)$] $|{\cal A}(x,z,\xi)-{\cal A}(y,w,\xi)|\leq \Lambda_1(1+g(|\xi|)(|x-y|^\alpha+|z-w|^\alpha)$
\item[$(d)$] $|B(x,z,\xi)|\leq \Lambda_1(1+g(|\xi|)|\xi|),$
\end{itemize}
for some positive constants $\Lambda$, $\Lambda_1$, $M_0$, for all $x$ and $y\in \Omega$, for all $z,w\in [-M_0,M_0]$ and for all $\xi \in \mathbb{R}^n$.	
Then, any solution $u\in W^{1,G}(\Omega)$, with $|u|\leq M_0$ in $\Omega$, is  $C^{1,\beta}(\overline{\Omega})$ for some positive $\beta$. Moreover 
\begin{align}\label{boundc1}
	\|u\|_{C^{1,\beta}(\overline{\Omega})}\leq C(\alpha, \Lambda, \delta, g_0, n, \Lambda_1, g(1), \Omega, M_0)
\end{align}
\end{proposition}

\begin{lemma}\label{LemLib} Let $A$ be a Young function satisfying $[A1]$ and $[A2]$. Put $\Phi(\xi)=A(|\xi|)$. Then 
\begin{eqnarray}\label{1.10a''}
		\sum_{i,j=1}^n\partial_{ij}\Phi(\eta)\xi_i\xi_j
		\geq\min\{\delta, 1\}\frac{A'(|\eta|)}{|\eta|}|\xi|^2\ \ \hbox{for all}\ \xi \in \rn,\eta\in \rn\setminus\{0\}\,,
\end{eqnarray}
\begin{eqnarray}\label{1.10b}
\sum_{i,j=1}^n\left|\partial_{ij}\Phi(\eta)\right|
%\leq\frac{\left(\sum_{i=1}^n|\eta_i|\right)^2}{{|\eta|^2}} \left|A''(|\eta|) %-\frac{A'(|\eta|)}{|\eta|}\right|+
%n\frac{A'(|\eta|)}{|\eta|}\nonumber \\
	\leq [2\max\{|\delta-1|,|g_0-1|\}+n]\frac{A'(|\eta|)}{|\eta|}\ \ \hbox{for all}\ \eta\in \rn\setminus\{0\}\,,
\end{eqnarray} 	
\end{lemma}
and $\nabla \Phi=\cal A$ satisfies conditions $a--c$ in  Proposition \ref{Lib}, with $g(t)=\min\{\delta,\,1\}A'(t)$, $\Lambda= \frac{\lambda}{\min\{\delta,\,1\}}$, where $\lambda=2\max\{|\delta-1|,|g_0-1|\}+n$.\\
%, and $(1.10)c$ holds whatever $\Lambda_1$ is.\\
{\bf Proof}
From \eqref{dg} 
\begin{equation*}
(\delta -1) \frac{A'(|\eta|)}{|\eta|}\leq  A''(|\eta|)-\frac{A'(|\eta|)}{|\eta|} \leq (g_0-1)\frac{A'(|\eta|)}{|\eta|}
\ \ \ \hbox{for all}\ \eta\in \rn\setminus\{0\}\,.
\end{equation*}
Also, $\partial_i\Phi(\eta)=A'(|\eta|)\frac{\eta_i}{|\eta|}$, and
\begin{equation*}
\partial_{ij}\Phi(\eta)= A''(|\eta|)\frac{\eta_i\eta_j}{|\eta|^2} +A'(|\eta|)\left(\frac{\delta_{ij}}{|\eta|} -\frac{\eta_i\eta_j}{|\eta|^3}\right)\ \ \hbox{for all}\ \eta\in \rn\setminus\{0\}\,.
\end{equation*}
Thus
\begin{align*}
\sum_{i,j=1}^n\partial_{ij}\Phi(\eta)\xi_i\xi_j= &\sum_{i,j=1}^n \left(A''(|\eta|)\frac{\xi_i\eta_i\xi_j\eta_j}{|\eta|^2} +A'(|\eta|)\frac{\delta_{ij}\xi_i\xi_j}{|\eta|}- A'(|\eta|)\frac{\xi_i\eta_i\xi_i\eta_j}{|\eta|^3}\right)\\
&=\left(\frac{A''(|\eta|)}{|\eta|^2}- \frac{A'(|\eta|)}{|\eta|^3}\right)(\xi,\eta)^2+
\frac{A'(|\eta|)}{|\eta|}|\xi|^2\\
 &\geq (\delta -1)\frac{A'(|\eta|)}{|\eta|^3}(\xi,\eta)^2+\frac{A'(|\eta|)}{|\eta|}|\xi|^2\ \ \ \hbox{for all}\ \xi \in \rn,\eta\in \rn\setminus\{0\}\,.\qquad\qquad\qquad\qquad
\end{align*}
If $\delta\geq 1$ 
\begin{equation}\label{1.10a}
    \sum_{i,j=1}^n\partial_{ij}\Phi(\eta)\xi_i\xi_j\geq \frac{A'(|\eta|)}{|\eta|}|\xi|^2\ \ \hbox{for all}\ \xi \in \rn,\eta\in \rn\setminus\{0\}\,.
\end{equation}
If $\delta<1$ 
\begin{eqnarray}\label{1.10a'}
    \sum_{i,j=1}^n\partial_{ij}\Phi(\eta)\xi_i\xi_j\geq (\delta -1)\frac{A'(|\eta|)}{|\eta|^3}|\xi|^2|\eta|^2+\frac{A'(|\eta|)}{|\eta|}|\xi|^2\\
    =\delta\frac{A'(|\eta|)}{|\eta|}|\xi|^2\ \ \hbox{for all}\ \xi \in \rn,\eta\in \rn\setminus\{0\}\,.\nonumber
\end{eqnarray}
Putting together \eqref{1.10a} and \eqref{1.10a'}, we get \eqref{1.10a''}.\\
Consider
\begin{equation*}
|\partial_{ij}\Phi(\eta)|\leq\frac{|\eta_i||\eta_j|}{{|\eta|^2}} \left|A''(|\eta|) -\frac{A'(|\eta|)}{|\eta|}\right|+
	\delta_{ij}\frac{A'(|\eta|)}{|\eta|} \quad \hbox{for all}\ \eta\in \rn\setminus\{0\}\,.
\end{equation*}	

Thus
\begin{eqnarray*}
\sum_{i,j=1}^n\left|\partial_{ij}\Phi(\eta)\right|\leq\frac{\left(\sum_{i=1}^n|\eta_i|\right)^2}{{|\eta|^2}} \left|A''(|\eta|) -\frac{A'(|\eta|)}{|\eta|}\right|+
n\frac{A'(|\eta|)}{|\eta|}\nonumber \\
\leq\left[2\max\{|\delta-1|,|g_0-1|\}+n\right]\frac{A'(|\eta|)}{|\eta|}\ \ \hbox{for all}\ \eta\in \rn\setminus\{0\}\,.
\end{eqnarray*} 
So \eqref{1.10b} holds with $\lambda=2\max\{|\delta-1|,|g_0-1|\}+n$.
%Conditions $(1.10)a--c$ in \cite{Li1} hold with 
%$a_{ij}=\partial_{ij}\Phi$, 
\qed
\section{Main results}\label{sec3}
In this section first we give two existence and regularity results (Theorems \ref{main1} and \ref{main2}), in which we assume a global growth condition on $f$, unilateral with respect to $s\in \r$. In Theorem \ref{main1} we require that $f$ satisfies some conditions for $(x,s,\xi)\in \Omega \times  [0,+\infty)\times \rn$ and obtain the existence of a nonnegative solution. Similarly, in Theorem \ref{main2}, $f$ satisfies some conditions for $(x,s,\xi)\in \Omega \times  (-\infty, 0]\times \rn$ that guarantee the existence of a non-positive solution.\\
\noindent Here is the definition of weak solution to \eqref{problem}. 
\begin{definition}
A function  $u\in \w0$ is a weak solution to problem \eqref{problem} if
\begin{equation*}
\int_{\Omega}A'(|\nabla u|)\cdot\frac{\nabla u}{|\nabla u|}\cdot\nabla v dx= \int_{\Omega}f(x,u,\nabla u)vdx
\end{equation*}
for all $v\in \w0$. 
\end{definition}
%We do not require $f(x,u,\nabla u)v\in L^1(\Omega)$ for all $u,v\in \w0$. From %(6.23) of \cite{BaCi1}, $A'(|\nabla u|)\cdot\frac{\nabla u}{|\nabla u|}\in %L^{\widetilde A}(\Omega)$ for all $u,v\in \w0$. Thus using \eqref{holderyoung} 
%$$\left|\int_{\Omega}A'(|\nabla u|)\cdot\frac{\nabla u}{|\nabla u|}\cdot \nabla %v dx\right|\leq 2\left\|A'(|\nabla u|)\frac{\nabla u}{|\nabla %u|}\right\|_{L^{\widetilde A}(\Omega}\|v\|_{L^{A}(\Omega)}\,.$$ 
%This guarantees that, if $u$ is a solution, then $f(x,u,\nabla u)v\in %L^1(\Omega)$ for every $v\in \w0$.\\
For the first two Theorems, we assume
$$({\cal{H}}):
\left\{
\begin{array}{l}
	a:[0,+\infty[\to [0,+\infty[ \ \hbox{is a locally essentially bounded function};\\
	\rho_1,\rho_2:\Omega\to [0,+\infty[\ \hbox{are two measurable functions},\ \rho_1,\,\rho_2\in L^\infty(\Omega)\ \hbox{and}\\ \rho_2 (x)>0\ \hbox{on a set of positive measure};\\
	g_1,g_2:[0,+\infty[\to [0,+\infty[\ \hbox{are two non-decreasing functions such that}\ g_1(0)=g_2(0)=0\\
	\hbox{and there exist}\ s_0>0,\,k_1\in \left ]0,\omega_n^{\frac{1}{n}}|\Omega|^{-\frac{1}{n}}\right[,\ \hbox{such that}\ 
	g_1(|s|)|s|\leq A(k_1|s|)\ \hbox{for all}\ |s|\geq s_0\,.
\end{array}
\right.
$$
Here $\omega_n$ is the measure of the unit ball in $\rn$.
\begin{theorem}\label{main1}
Let $\Omega$ be a bounded domain in $\rn$ with $C^{1,\alpha}$ boundary. 
%for some $\alpha>0$. 
Let $A: [0,+\infty[\to [0,+\infty[$ be a Young function, satisfying $[A1]$ and $[A2]$. Let $f:\Omega\times \mathbb{R}\times \mathbb{R}^n\to \mathbb{R}$ be a Carath\'{e}odory function fulfilling
	\begin{equation}\label{growth f2+}
		\rho_2(x)-g_2(s)-a(s)A'(|\xi|)|\xi|\leq f(x,s,\xi)\leq \rho_1(x)+g_1(s)\ \hbox{ for a.e.}\ x\in \Omega,\ \hbox{all}\ s\geq 0,\ \hbox{all}\ \xi\in\rn\,.
		%(s,\xi)\in \r\times\rn\,.
	\end{equation}
The functions $a,\,\rho_1,\,\rho_2,\,g_1,\,g_2$ are as in $(\cal{H})$.
	Then problem \eqref{problem} has a nontrivial, nonnegative solution $u\in C_0^{1,\beta}(\overline \Omega)$.\\
	%for some $\beta>0$
	If, in addition, there exist $\overline{\delta} >0$ and $k_3>0$ such that
	$g_2(s)s\leq A(k_3s)$ for every $s\in (0,\overline{\delta})$, then $u>0$ in $\Omega$.
\end{theorem}
\begin{theorem}\label{main2}
Let $\Omega$ be a bounded domain in $\rn$ with $C^{1,\alpha}$ boundary.
%for some $\alpha>0$. 
Let $A: [0,+\infty[\to [0,+\infty[$ be a Young function, satisfying $[A1]$ and $[A2]$. Let $f:\Omega\times \mathbb{R}\times \mathbb{R}^n\to \mathbb{R}$ be a Carath\'{e}odory function fulfilling
\begin{equation}\label{growth f2}
-\rho_1(x)-g_1(|s|)	\leq f(x,s,\xi)\leq -\rho_2(x)+g_2(|s|)+a(s)A'(|\xi|)|\xi|\ \hbox{for a.e.}\, x\in \Omega,\ \hbox{all}\ s\leq 0,\ \hbox{all}\ \xi\in\rn\,,
	%(s,\xi)\in \r\times\rn\,.
	\end{equation}
where the functions $a,\,\rho_1,\,\rho_2,\,g_1,\,g_2$ are as in $(\cal{H})$.
	Then problem \eqref{problem} has a nontrivial, non-positive solution $u\in C_0^{1,\beta}(\overline \Omega)$.\\
	%for some $\beta>0$.\\
	If, in addition, there exist $\overline{\delta} >0$ and $k_3>0$ such that
	$g_2(s)s\leq A(k_3s)$ for every $s\in (0,\overline{\delta})$, then $u<0$ in $\Omega$.
\end{theorem}
\begin{remark}
In \cite{BalFil}, Theorem 1, the authors prove the existence of a positive solution for a problem with the $p$-Laplacian, and a convection term $f$ satisfying the hypotheses of Theorem \ref{main2}.
\end{remark}
%$\rho_1=\rho_2=g_1\equiv 0, g_2(t)=(t^+)^q+Ct^s$ e l'esponente nel gradiente è minore di p.
For the proof of the Theorems above, we need an abstract existence result, where sub and supersolutions come into play.\\
The definition of sub and supersolution in general domains, for which a trace theory may not hold, can be found in \cite{BaTo1}. Our hypotheses on $\Omega$ allow to adopt the classical definition (see \cite[Theorem 3.1]{Ctrace}).\\
%onsider sub and supersolutions belonging in $W^{1,\infty}(\Omega)$, where
We say that $\overline u\in W^{1,A}(\Omega)$ is a supersolution to \eqref{problem} if $\overline u_{|\partial \Omega}\geq 0$ (in the sense of traces) and
\begin{equation*}\label{supersolution}
	\int_{\Omega}A'(|\nabla\overline u|)\cdot\frac{\nabla\overline u}{|\nabla\overline u|}\cdot\nabla v dx\geq \int_{\Omega}f(x,\overline u,\nabla \overline u)vdx
\end{equation*}
for all $v\in W_0^{1,A}(\Omega)$, $v\geq 0$ a.e. in $\Omega$.\\
We say that $\underline u\in W^{1,A}(\Omega)$ is a subsolution to \eqref{problem} if $\underline u_{|\partial \Omega}\leq 0$ (in the sense of traces) and
\begin{equation*}\label{subsolution}
	\int_{\Omega}A'(|\nabla\underline u|)\cdot\frac{\nabla \underline u}{|\nabla \underline u|}\cdot \nabla v dx\leq \int_{\Omega}f(x,\underline u,\nabla \underline u)vdx
\end{equation*}
for all $v\in W_0^{1,A}(\Omega)$, $v\geq 0$ a.e. in $\Omega$.\\
For the next Theorem we assume that problem \eqref{problem} has a subsolution and supersolution, $\underline{u}$, $\overline{u}\in W^{1,\infty}(\Omega)$, with $\underline{u}(x)<\overline{u}(x)$ for all $x\in \Omega$.
Also, $f:\Omega\times \r\times \rn\to \r$ is a Carath\'{e}odory function satisfying the following growth condition:
\begin{itemize}
\item [(H)] there exists a function $\sigma\in L^{\infty}(\Omega)$ and a constant $a>0$, such that
$$ %\label{growthf}
|f(x,s,\xi)|\leq \sigma (x)+a A'(|\xi|)|\xi|\quad \hbox{for a.e.}\ x\in \Omega,\ \hbox{all}\ s\in [\underline u(x),\overline u(x)],\ \hbox{all}\ \xi\in \rn\, .
$$
\end{itemize}
\noindent The local condition on $f$, with respect to $s$, is sufficient for our purposes. The use of the method of sub and supersolutions requires an a priori analysis of the problem. Only once the existence of sub and supersolutions has been established does one proceed to search for the existence of a solution.

\begin{theorem}\label{main}
Let $\Omega$ be a bounded domain in $\rn$ with $C^{1,\alpha}$ boundary. %$\partial \Omega$ for some $\alpha>0$. 
Let $A: [0,+\infty[\to [0,+\infty[$ be a Young function satisfying $[A1]$ and $[A2]$. Let $\underline{u}$, $\overline{u}\in W^{1,\infty}(\Omega)$ be as above and assume that $f$ satisfies hypothesis (H). Then problem \eqref{problem} admits at least a solution  $u\in C_0^{1,\beta}(\overline \Omega)$. Moreover $\underline u(x)\leq u(x)\leq \overline u(x)$ a.e in $\Omega$.
\end{theorem}
{\bf Proof}. Let $M=\max\{\| \overline u\|_\infty,\| \underline u\|_\infty\}$ and $R>\max\{\|\nabla \overline u\|_\infty,\|\nabla \underline u\|_\infty\}$. Consider the truncated function $f_R$ defined by
\begin{equation*}\label{fR}
    f_R(x,s,\xi)=\left\{\begin{array}{cc}
                          f(x,s,\xi) & \hbox{if}\ |\xi|\leq R,\\
                          f(x,s,\xi)\cdot\frac{A'(R)R}{A'(|\xi|)|\xi|} & \hbox{if}\ |\xi|> R \,,
                        \end{array}
    \right.
\end{equation*}
and the problem
\begin{equation}\label{problem1}
\begin{cases}
%- {\rm div}(A'(|\nabla u|)\frac{\nabla u}{|\nabla u|}) 
-\Delta_A(u)=f_R(x,u,\nabla u) & {\rm in}\,\,\,\Omega \\
u  =0 & {\rm on}\,\,\,\partial \Omega \,.
\end{cases}
\end{equation}
In view of the choice of $R$, $\underline u$ and $\overline u$ are a subsolution and a supersolution to \eqref{problem1} respectively. Using the monotonicity of $A'$ we deduce that $|f_R(x,s,\xi)|\leq \sigma (x)+a A'(R)R$, for a.e. $x\in \Omega$, all $s\in [\underline u(x),\overline u(x)]$, all $\xi\in \rn$. From Theorem 3.6 in \cite{BaTo1} problem \eqref{problem1} admits a solution $u\in W^{1,A}_0(\Omega)$ with $\underline u(x)\leq u(x)\leq \overline u(x)$ a.e in $\Omega$. Thus $u\in L^\infty(\Omega)$.\\ 
\noindent The functions $A$ and $f$ satisfy the hypotheses of Proposition \ref{Lib}, with $\Lambda_1=\max\{\|\sigma\|_\infty, \min\{\delta,1\}^{-1}a\}$ (see also Lemma \ref{LemLib}). Since $|f_R|\leq |f|$ the same holds for $f_R$ whatever $R$ is.
Due to Proposition \ref{Lib} there exist two positive constants $0<\beta\leq 1$ and $C$,  independent from $R$, such that any solution to \eqref{problem1} belongs in $C_0^{1,\beta}(\overline \Omega)$ and $\|u\|_{C_0^{1,\beta}(\overline \Omega)}\leq C$. Choosing $R>C$ we deduce that $u$ is a solution to \eqref{problem}.
\qed
%\begin{remark}
%The proof above shows also that under the weaker assumption %$\sigma\in L^{\widetilde A_n}(\Omega) $, the solution $u\in %L^\infty(\Omega)$. This is due to the fact that the upper bound %$|f_R(x,s,\xi)|\leq \sigma (x)+a A'(R)R$, for a.e. $x\in \Omega$, %all $s\in [\underline u(x),\overline u(x)]$, all $\xi\in \rn$, %still allows to apply Theorem 2.11 in \cite{BaTo1}.
%\end{remark}
\begin{remark}
When the solution $u$ has constant sign, it should be interest to verify if it is positive (or negative) in $\Omega$. The maximum principle by Pucci-Serrin (see \cite [Theorem 3.5] {PS}) is a powerful tool, as it ensures that (under some proper conditions on $A$ and $f$) 
%When $f$ satisfies:
%\begin{equation}\label{growthPS}
%	f(x,s,\xi)\geq -a A'(|\xi|)-b(s)\ \hbox{for all} \; x\in \Omega,\ s>0,\ %\xi\in \mathbb{R}^n,\; |\xi|\leq 1,
%\end{equation}
%and 
%\begin{equation}\label{divPS}
%	\int_0\frac{1}{H^{-1}(B(s))}ds=+\infty\,,
%\end{equation}
%where $b:[0,+\infty[ \to [0,+\infty[$ is a function increasing in %$(0,\overline{\delta})$ (for some $\overline{\delta}>0$), and $b(0)=0$, 
any nonnegative solution to \eqref{problem} is positive in $\Omega$. 
%Here, as usual, $B(s)=\int_0^s b(\tau)d\tau$, while $H(s)=s A'(s)-A(s)$.
A quite standard situation occurs when $f$ is bounded below by suitable monotone functions and $A\in\Delta_2$ near $0$, as the following corollary shows. \end{remark}
\begin{corollary}\label{cor1}
Under the hypotheses of Theorem \ref{main}, assume that $f$ satisfies
\begin{equation}\label{growthPS}
	f(x,s,\xi)\geq -a A'(|\xi|)-b(s)\ \hbox{for all} \; x\in \Omega,\ s>0,\ \xi\in \mathbb{R}^n,\; |\xi|\leq 1,
\end{equation}
where $a>0$, $b:[0,+\infty[ \to [0,+\infty[$ is a function increasing in $(0,\overline{\delta})$ (for some $\overline{\delta}>0$), $b(0)=0$, and $b(s)=\frac{A(ks)}{s}$ for $s\in (0,\overline\delta)$ and some $k>0$.
Then any nonnegative, nontrivial solution to \eqref{problem} is positive.
\end{corollary}
{\bf Proof.}
Let $u\in W^{1,A}_0(\Omega)$ be a nonnegative, nontrivial solution to \eqref{problem}. Theorem \ref{main} ensures that $u\in C_0^{1,\beta}(\overline{\Omega})$. In order to prove that $u>0$ in $\Omega$ we use Theorem 5.3.1 of \cite{PS}.\\
Conditions $(A1)'$ and $(A2)$ of the Theorem cited above hold, because $A\in C^2((0,+\infty))$, $A'(0)=0$, and $s\mapsto A'(s)$ is strictly increasing. Conditions $(F2)$ and $(B1)$ are satisfied too.\\
It remains to verify condition $(1.1.5)$ of Theorem 5.3.1 of \cite{PS}. 
Put $B(s)=\int_0^s b(t)dt$. Due to the monotonicity of $\frac{A(t)}{t}$, for $s\in(0,\overline\delta)$, it holds 
$$B(s)=\int_0^s \frac{A(kt)}{t}dt\leq \int_0^s \frac{A(ks)}{s}dt=A(ks)\,.$$
% and  $sA'(s)\in C^1(]0,+\infty[)$ Also,the function $b(s)=\frac{A(k_3s)}{s}$ for all %$s>0$ it is increasing and $\lim_{s\to 0^+}b(s)=0$. The function $kb(s)$ satisfies %$(F2)$. From \eqref{growth f2} and taking into account additional condition on $g_2$ %we obtain
If $h\in \mathbb{N}$ is such that $k<2^h$ then, in view of \eqref{delta2young}
$$A(ks)\leq K^h A(s)\quad \hbox{for all}\; s\geq 0\,.$$ 
Let $b_1=\max\{p_A-1,K^h\}$. 
%$$
%f(x,s,\xi)\geq \rho_2(x)-g_2(|s|)-a A'(|\xi|)||\xi|\geq -b(|s|)-a A'(|\xi|)$$ for all %$x\in \Omega$, all $s\in \mathbb{R}$, and all $\xi\in \mathbb{R}^n$, $|\xi|\leq 1$. 
Then, for $s\in (0,\overline{\delta})$, using \eqref{convexity} and the inequality above
%,  $\frac{p_A-1}{b_1}\leq 1$
\begin{align*}
	H(s)=sA'(s)-A(s)&\geq (p_A-1)A(s)=\frac{p_A-1}{b_1}b_1A(s)\nonumber\\
	\geq  b_1 A\left(\frac{(p_A-1)s}{b_1}\right)\geq  K^hA\left(\frac{(p_A-1)s}{b_1}\right)&\geq  A\left(\frac{k(p_A-1)s}{b_1}\right) \geq B\left(\frac{(p_A-1)s}{b_1}\right)\,, 
\end{align*}
or equivalently
\begin{equation*}
H\left(\frac{b_1s}{p_A-1}\right)\geq B(s)\quad \hbox{for}\ 0<s<\frac{b_1\overline{\delta}}{p_A-1}\,.
\end{equation*}
$H$ is increasing, so
\begin{equation*}
	\frac{b_1s}{p_A-1}\geq H^{-1}(B(s))\quad \hbox{for}\ 0<s<\frac{b_1\overline{\delta}}{p_A-1}\,.
\end{equation*}
Finally 
\begin{equation}\label{minor per strong}
\frac{1}{ H^{-1}(B(s))}\geq \frac{p_A-1}{b_1s}\quad \hbox{for}\ 0<s<\frac{b_1\overline{\delta}}{p_A-1}\,.
\end{equation}
Integrating \eqref{minor per strong} from $\ep$ to $s<\frac{b_1\overline{\delta}}{p_A-1}$ and passing to the limit as $\ep\to 0^+$ we obtain condition $(1.1.5)$ of Theorem 5.3.1 of \cite{PS}. Thus $u>0$ in $\Omega$.\qed

Now we accomplish with the proof of Theorems \ref{main1} and \ref{main2}.\\
{\bf Proof of Theorem \ref{main1}}. 
From the proof of Theorem 3.3 of \cite{BaTo1} we know that there exists a nontrivial solution $\overline{u}\geq 0$, to the problem
	$$
	- \Delta_A(u) =\rho_1(x)+g_1(|u|)\,.
	$$
Theorem 3 of \cite{Cbound} guarantees that $\overline{u}$ is bounded. Finally, from Proposition \ref{Lib}, we have that $\overline{u}\in C_0^{1,\beta}(\overline{\Omega})$.
The inequalities in \eqref{growth f2+} show that $\overline{u}$ is a supersolution to problem \eqref{problem} and $\underline{u}=0$ is a subsolution to problem \eqref{problem}. The assumptions on $\rho_2$ guarantee that $\underline{u}=0$ is not a solution. If we put $\sigma(x)=\max\{\rho_1(x)+g_1(\overline{u}(x)) ,\; \rho_2(x)+g_2(\overline{u}(x))\}$ for a.e. $x\in \Omega$, then \eqref{growth f2+} leads to
$$
|f(x,s,\xi)|\leq \sigma(x)+a(s)A'(|\xi|)|\xi|\quad \hbox{for a.e.} \; x\in \Omega, \; \hbox{all}\ s\in [0,\overline{u}(x)],\; \hbox{all}\ \xi\in \mathbb{R}^n\,.$$
Let $I=[0,\sup_\Omega \overline u]$ and $a=\|a\|_{L^\infty(I)}$. Then $a<+\infty$ and $a(s)\leq a$ for a.e. $s\in I$. Let $I_0\subset I$ be a set of null measure, such that $a(s)>a$ for all $s\in I_0$. For $s\in I_0$ it holds
\begin{align}
|f(x,s,\xi)|=\lim_{t\to s}|f(x,t,\xi)|=\liminf_{t\to s}|f(x,t,\xi)|& \leq \sigma(x)+\liminf_{t\to s}a(t)A'(|\xi|)|\xi|\nonumber\\
	\leq \sigma(x)+aA'(|\xi|)|\xi| &\quad \hbox{for a.e.} \; x\in \Omega, \; \xi\in \mathbb{R}^n\,.
\end{align}
Thus  
\begin{align}
|f(x,s,\xi)|\leq \sigma(x)+aA'(|\xi|)|\xi|&\quad \hbox{for a.e.} \; x\in \Omega, \; s\in [0,\overline u(x)],\; \xi\in \mathbb{R}^n\,.
\end{align}
So, from Theorem \ref{main}, problem \eqref{problem} admits at least a nontrivial solution $u\in C_0^{1,\beta}(\overline{\Omega})$ such that $0\leq u\leq \overline{u}$.\\
Now we prove that, under the additional condition on $g_2$, $u>0$ in $\Omega$. We set $b(s)=\max \{g_2(s),\frac{A(k_3s)}{s}\}$, for $s\geq 0$ and observe that the left inequality in \eqref{growth f2+} guarantees that we can apply Corollary \ref{cor1}.
\qed
\noindent {\bf Proof of Theorem \ref{main2}}. It is enough to put $f_1(x,s,\xi)=-f(x,-s,-\xi)$ and to use Theorem \ref{main1} for $f_1$.\qed
% L'operatore $A_P$ e' il nostro $\frac{A'(s)}{s}$, la $f_p$ e' la nostra $kb$, la $F$  .

\section{Examples}\label{sec4}
This section is devoted to some examples with different Young functions and various nonlinearities.\\
% The functions $A$ involved fall into the more general class of those functions %$A(t)\approx t^p\lg^\alpha(c+t)$ near infinity, for $p>1$ and suitable $\alpha, %c\in \r$. 
\noindent In the first two examples we consider the problem
\begin{equation}\label{pqf}
\ \begin{cases}
- {\rm div}\left((|\nabla u|^{p-2}+|\nabla u|^{q-2})\nabla u\right) =f(x,u,\nabla u) & {\rm in}\,\,\, \Omega \\
u  =0 & {\rm on}\,\,\,\partial \Omega \, ,
\end{cases}
\end{equation}
where $1<q<p<\infty$. Here $A(t)=\frac{t^p}{p}+\frac{t^q}{q}$ for all $t\geq 0$, and 
%q>1 serve per delta 0>0 a continua serve per carateodory
\eqref{dg} holds with $\delta=q-1$, $g_0=p-1$.
\begin{example}
Let $a:\r\to\r$ and $g:\r\to\r$ be two continuous functions and let $h:\Omega\to \r$ be an essentially bounded function.
Assume that there exist $s_1,\, s_2\in\r$, such that $g(s_1)=g(s_2)=0$ for some $s_1<0<s_2$, $g(s)\neq 0$ for all $s\in ]s_1,s_2[$, and $|\{x\in\Omega:h>0\}|>0$, $|\{x\in\Omega:h<0\}|>0$.\\
Set
\begin{equation*}
f(x,s,\xi)= g(s)h(x)+a(s)A'(|\xi|)(|\xi|)\ \hbox{for}\ (x,s,\xi)\in \Omega \times \r\times \rn\,.
\end{equation*}
Then $u_1=s_1$ and $u_2=s_2$ are a subsolution and a supersolution to \eqref{pqf}, respectively. Also $u\equiv 0$ is not a solution nor a sub or a supersolution and $f$ satisfies condition (H) 
 with $\sigma (x)=|h(x)|\max_{[s_1,s_2]}|g(s)|$ and $a =\max_{[s_1,s_2]}|a(s)|$. By Theorem \ref{main}, problem \eqref{pqf} has a nontrivial solution $u\in C^{1,\beta}_0(\overline\Omega)$ with $s_1\leq u\leq  s_2$ a.e. in $\Omega$.

\end{example}

\begin{example}
Let $a:\r\to\r$ and $g:\r\to\r$ be two continuous functions and let $h:\Omega\to \r$ be an essentially bounded function. Assume that $h(x)\geq 0$ (or $h(x)\leq 0$) in $\Omega$, $h(x)>0$ on a set of positive measure, $g(s_1)=0$ for some $s_1>0$, $g(s)\neq 0$ for all $s\in [0,s_1[$, and $g(0)h(x)\geq 0$ in $\Omega$.\\
Set
\begin{equation*}
	f(x,s,\xi)= g(s)h(x)+a(s)A'(|\xi|)(|\xi|)\ \hbox{for}\ (x,s,\xi)\in \Omega \times \r\times \rn\,.
\end{equation*}
Then $u_1\equiv0$ and $u_2\equiv s_1$ are a subsolution and a supersolution to \eqref{pqf}, respectively. Also $u\equiv 0$ is not a solution and $f$ satisfies condition (H) with $\sigma (x)=|h(x)|\max_{[0,s_1]}|g(s)|$ and $a =\max_{[0,s_1]}|a(s)|$. By Theorem \ref{main}, problem \eqref{pqf} has a nontrivial solution $u\in [0,s_1]$, $u\in C^{1,\beta}_0(\overline\Omega)$. From Theorem 5.3.1 in \cite{PS}, $u>0$ in $\Omega$ (note that $g(s)h(x)\geq 0$ in $\Omega\times[0,s_1]$).
%b(s)=0 perche g(s)h(x)\geq 0 in (o,\delta)
\end{example}
\begin{example} 
Consider the problem 
\begin{equation}\label{problempq}
\begin{cases}
- {\rm div}(\lg(1+|\nabla u|^{q})|\nabla u|^{p-2}\nabla u) =f(x,u,\nabla u) & {\rm in}\,\,\, \Omega \\
u  =0 & {\rm on}\,\,\,
\partial \Omega \,,
\end{cases}
\end{equation}
with $q>0$, $p>1$.\\ 
Let $r\geq q$, $\overline \delta >0$ and define $b:[0,+\infty[\to [0,+\infty[$ as
\begin{equation*}
    b(s)=\left\{\begin{array}{cc}
         s^{p+q-1}&\hbox{if}\ s\in [0,\overline \delta],  \\
          s^{p+r-1}& \hbox{if}\ s> \overline \delta\,.
    \end{array}\right.
\end{equation*}
Assume that $f:\Omega\times\r\times\rn\to\r$ is a Carath\'{e}odory function satisfying
\begin{itemize}
\item [$(f_0)$] there exists $\sigma >0$ such that $f(x,\sigma,0)\leq 0$ a.e. in $\Omega$;
\item [$(f_1)$] $f(x,0,0)\geq 0$ a.e. in $\Omega$, with strict inequality on a set of positive measure;
\item [$(f_2)$] there exist $k>0$ such that
 \begin{equation*}\label{strong1}
f(x,s,\xi)\geq-k(|\xi|^{p-1}\lg(1+|\xi|^q)+b(s))\ \hbox{for}\ x\in\Omega,\ \hbox{all}\ s\geq 0\ \hbox{and all}\ \xi\in \rn,\ \hbox{with}\ |\xi|\leq 1\,.
  \end{equation*}
\item [$(f_3)$]  there exists $c >0$ such that 
$|f(x,s,\xi)|\leq c(1+|\xi|^p\lg(1+|\xi|^q))\ \hbox{for all}\ x\in\Omega,\ s\in \r, \ \xi\in\rn$.
\end{itemize}
For the function $A'$ in problem \eqref{problempq} condition \eqref{dg} holds with $\delta=p-1$, $g_0=p-1+q$ and 
%$A'(t)=t^{p-1}\lg(1+t^{q})$ and 
$A(t)\approx  t^{p+q}$ for  $t$ small.
%for some $\overline{\delta}>0$. 
We can apply Theorem \ref{main} and Corollary \ref{cor1} to obtain the existence of a positive solution $u\in C_0^{1,\beta}(\overline \Omega)$, and $u\leq \sigma$ in $\Omega$.
\end{example}
The example above extends in two directions Theorem 6 of \cite{FMMT}: it allows an higher growth for $f$ and permit also the choice $r=q$ in the lower bound for $f$.

\begin{example} 
Consider the problem \begin{equation}\label{Apq}
	\begin{cases}
		- {\rm div}\left(\frac{|\nabla u|^{p-2}\nabla u}{\lg^q(1+|\nabla u|)}\right) =f(x,u,\nabla u) & {\rm in}\,\,\, \Omega \\
		u  =0 & {\rm on}\,\,\,
		\partial \Omega \,,
	\end{cases}
\end{equation}
with $p>1$, $p-q-1>0$.
% delta =p-1-q, g_0=p-1
Let $\rho\in L^{\infty}(\Omega)$ and $g_1,g_2:[0,+\infty[\to [0,+\infty[$ be two unbounded, nondecreasing functions, such that $g_1(0)=g_2(0)=0$. Also, let $a_1,a_2: \r \to [0,+\infty[$ be two locally essentially bounded functions and let $c_1,c_2>0$.\\
Assume that  $f:\Omega\times\r\times\rn\to\r$ is a Carath\'{e}odory function satisfying
\begin{eqnarray*}
-c_1+g_1(|s|)-a_1(s)\frac{|\xi|^{p}}{\lg^{q}(1+|\xi|)}\leq f(x,s,\xi)\leq -c_2+g_2(|s|)\rho(x)+a_2(s)\frac{|\xi|^{p}}{\lg^{q}(1+|\xi|)}\
\end{eqnarray*}
$\hbox{for}\ (x,s,\xi)\in \Omega \times \r\times \rn \,.$ \\
We show that problem \eqref{Apq} has a nontrivial solution $u\leq 0$ in $\Omega$.\\
For the function $A'$ in problem \eqref{Apq} condition \eqref{dg} holds with $\delta =p-1-q$, $g_0=p-1$.
If $k:=\inf\{s>0\,:\, g_1(s)\geq c_1\}$, then $\underline u\equiv -k$ is a subsolution to \eqref{Apq}, and $\overline u\equiv 0$ is a supersolution but not a solution to \eqref{Apq}. Let $a=max \{\|a_1\|_{L^\infty ([-k,0])},\ \|a_2\|_{L^\infty ([-k,0])}\}$, $\sigma(x)=max\{c_1, -c_2+g_2(k)\rho(x)\}$. Then 
\begin{equation*}
|f(x,s,\xi)|\leq \sigma(x)|+a A'(|\xi|)|\xi|\quad \hbox{for}\ x\in \Omega,\ s\in\,[-k,0],\ \xi\in\rn\,.
\end{equation*}
%and $f$ satisfies $(H)$ with $\sigma (x)=b|\rho(x)|+c$. 
By Theorem \ref{main}, problem \eqref{Apq} has a nontrivial solution $u\in [-k,0]$.
\end{example}


\begin{thebibliography}{99}

\bibitem[Ad]{Adams} R.A. Adams, Sobolev Spaces, \emph{Academic Press}, New York, 1975.
\bibitem[BF]{BalFil} L. Baldelli, R. Filippucci, \emph{Existence results for elliptic problems with gradient terms via a priori estimates}, Nonlinear Anal. {\bf 198} (2020), 111894, 22 pp. 
\bibitem[BaCi1]{BaCi1} G. Barletta, A. Cianchi, \emph{Dirichlet problems for fully anisotropic elliptic equations}, Proc. Roy. Soc. Edinburgh Sect. A {\bf 147} 1 (2017), 25--60.
\bibitem[BaTo1]{BaTo1} G. Barletta, E. Tornatore, \emph{Elliptic problems with convection terms in Orlicz spaces}, J. Math. Anal. Appl. {\bf 495} (2021), 124779.
 \bibitem[BNV]{BNV} M.F. Bidaut-Véron, Marie-Françoise, Q.H. Nguyen, L. Véron, \emph{Quasilinear elliptic equations with a source reaction term involving the function and its gradient and measure data}, Calc. Var. Partial Differential Equations {\bf 59} 5, 148 (2020), 38 pp. 
\bibitem[Chl]{Chl} I. Chlebicka, \emph{A pocket guide to nonlinear differential equations in Musielak-Orlicz space}, Nonlinear Analisys {\bf 175} (2018), 1--27.
%\bibitem[Ci1]{Csharp} A. Cianchi, \emph{A sharp embedding theorem for Orlicz-Sobolev spaces}, Indiana Univ. Math. J. {\bf 45} (1996), 39--65.
\bibitem[Ci2]{Cbound} A. Cianchi, \emph{Boundedness of solutions to variational problems under general growth conditions}, Comm. Part. Diff. Eq. {\bf 22} (1997), 1629--1646.
%\bibitem[Ci3]{Cfully} A. Cianchi, \emph{A fully anisotropic Sobolev inequality}, Pacific J. Math. {\bf 196} (2000), 283--295.
\bibitem[Ci4]{Ctrace} A. Cianchi, \emph{Orlicz-Sobolev boundary trace embeddings}, Math. Z., {\bf 266} 2 (2010),  431--449.
\bibitem[FMP]{FMP} F.Faraci, D.Motreanu, D.Puglisi, Positive solutions of quasi-linear elliptic equations with dependence on the gradient, \emph{Calc. Var. Partial Differential Equations} {\bf 54} 1 (2015), 525--538.
\bibitem[FMMT]{FMMT} L.F.O. Faria, O.H. Miyagaki, D. Motreanu and M. Tanaka, \emph{Existence results for nonlinear elliptic equations with Leray-Lions operator and dependence on the gradient}, Nonlinear Anal. {\bf 96} (2014), 154--166.
\bibitem[FM]{FM} G.M. Figueiredo, G.F. Madeira, \emph{Positive maximal and minimal solutions for non-homogeneous elliptic equations depending on the gradient}, J. Differential Equations {\bf 274} (2021), 857–875. 
\bibitem[GW]{GW}  L.Gasi\'{n}ski, P.Winkert, Existence and uniqueness results for double phase problems with
convection term, \emph{J. Differential Equations}, https://doi.org/10.1016/j.jde.2019.10.022
\bibitem[G]{G} N.Grenon, Existence and Comparison Results for Quasilinear Elliptic Equations with Critical Growth in the Gradient, \emph{J. Differential Equations} {\bf 171} 1 (2001), 1--23.
%\bibitem[HLZ]{HLZ}  P. Hajlasz, Z. Liu, Zhuomin, \emph{A compact embedding of a Sobolev space is equivalent to an embedding into a better space},  Proc. Amer. Math. Soc.,  {\bf 138} 9 (2010), 3257--3266.
\bibitem[KrRu]{KrRu} M.A.Krasnosel'skii, Ja.B.Rutickii, "\emph{ Convex functions and Orlicz spaces}, Groningen: Noordhoff, 1961.
\bibitem[Li]{Li}  G.M. Lieberman, \emph{Boundary regularity for solutions of degenerate elliptic equations.} Nonlinear Anal. {\bf 12} (1988), 1203--1219.
\bibitem[Li1]{Li1} G.M. Lieberman, \emph{The natural generalization of the natural conditions of Ladyzhenskaya and Ural'tseva for elliptic equations},  Comm. Partial Differential Equations {\bf 16} (1991), 311--361.
\bibitem[MW]{MotWin} D. Motreanu, P. Winkert, \emph{Existence and asymptotic properties for quasilinear elliptic equations with gradient dependence}, Appl. Math. Lett. {\bf 95} (2019), 78–84.
\bibitem[NS]{NgSch} L.H.Nguyen, K. Schmitt, \emph{Boundary value problems for singular elliptic equations}, Rocky Mountain J. Math. {\bf 41} (2011), no. 2, 555–572.
\bibitem[PS]{PS} P. Pucci, J. Serrin, The maximum principle, Birkh\"{a}user Verlag, Basel 2007.
\bibitem[R]{Ruiz} D. Ruiz, \emph{A priori estimates and existence of positive solutions for strongly nonlinear problems}, J. Differential Equations {\bf 199} (2004), no. 1, 96–114. 
\bibitem[RR1]{RR1} M.M.Rao \& Z.D.Ren, ``Theory of Orlicz spaces", Marcel Dekker, New York, 1991.
\bibitem[RR2]{RR2} M.M.Rao \& Z.D.Ren, ``Applications of Orlicz spaces", Marcel Dekker, New York, 2002.
\bibitem[T]{Tanaka} M. Tanaka, \emph{Existence of a positive solution for quasilinear elliptic equations with nonlinearity including the gradient}, Bound. Value Probl.  {\bf  2013} 173 (2013), 11 pp. 
\bibitem[Z]{Zou} H.H. Zou, \emph{A priori estimates and existence for quasi-linear elliptic equations}, Calc. Var. Partial Differential Equations  {\bf 33} 4 (2008), 417–437.


\end{thebibliography}
\end{document}